\documentclass[11pt]{article}
\usepackage{latexsym}
\usepackage{amsfonts,amssymb,amsmath}
\newtheorem{thm}{Theorem}[section]
\newtheorem{cor}[thm]{Corollary}
\newtheorem{lem}[thm]{Lemma}
\newtheorem{pro}[thm]{Proposition}
\newtheorem{defn}[thm]{Definition}

\title{ Monoids that map onto the Thompson-Higman groups }

\author{ Jean-Camille Birget  }

\date{\today}
\begin{document}
\maketitle

\begin{abstract}
A slight modification of the definition of the Thompson-Higman groups
$G_{k,1}$ and $F_{k,1}$ leads to inverse monoids that map onto $G_{k,1}$
(respectively $F_{k,1}$), and that have interesting properties: they are
finitely generated, and residually finite.  These inverse monoids are 
closely related to the suffix expansion of $G_{k,1}$ (respectively 
$F_{k,1}$).  
\end{abstract}


\section{Introduction}

The groups $G_{k,1}$ and $F_{k,1}$ $(k \geq 2)$ of Thompson and Higman 
\cite{Th0,McKTh,Th}, \cite{Hig74} are well known in combinatorial group 
theory (see e.g.\ the references in \cite{BiThomps,BiThomMon}). 
A classical survey is \cite{CFP}. The group $G_{2,1}$ is usually called 
$V$, and $F_{2,1}$ is called $F$.

The group $G_{k,1}$ can be defined as consisting of all maximally 
extended right-ideal isomorphisms between finitely generated essential 
right ideals of a free monoid $A^*$ (where $A$ is a $k$-letter alphabet).
An important fact about isomorphisms between essential 
right ideals of $A^*$ is that they have a {\em unique} maximal essential
extension. Multiplication in $G_{k,1}$ is composition, followed by 
maximum extension. Details appear in the Background definitions below. 
In \cite{CFP} $G_{k,1}$ is defined (differently, but isomorphically) as 
consisting of all continuous, increasing, piecewise-linear bijections 
(with finitely many pieces), from the real interval $[0,1]$ onto itself, 
where only base-$k$ rationals are allowed as coordinates of articulation 
points. A base-$k$ rational is a rational number of the form $a/k^n$ 
with $a, n$ integers and $n \geq 0$. In this definition of $G_{k,1}$, 
maximum extension is achieved by repeatedly combining certain 
adjacent linear pieces into one linear piece if these pieces have 
the same slope, and if the domain intervals of these pieces are of 
the form 
 \ $[ \frac{ak}{k^{n+1}}, \frac{ak+1}{k^{n+1}}]$\,, 
 \ $[\frac{ak+1}{k^{n+1}}, \frac{ak+2}{k^{n+1}}]$\,, \ \ $\ldots$ , 
 \ $[\frac{ak+k-1}{k^{n+1}}, \frac{ak+k}{k^{n+1}}]$\,; 
 \ then $[\frac{a}{k^n}, \frac{a+1}{k^n}]$  is the domain interval 
of the corresponding combined piece. 

In order to gain a better understanding of the role of the maximum
essential extension one can look at what happens if maximum essential 
extension is simply left out in the definition of the Thompson-Higman
groups.  The structure that is obtained then is the 
inverse monoid ${\sf riAut}(k)$, consisting of all isomorphisms 
between finitely generated essential right ideals of $A^*$, where 
$|A| = k$.  Multiplication is now just composition.
In relation to the group $F_{k,1}$ one can also define the inverse 
monoid ${\sf riAut}_{\sf dict}(k)$, consisting of all dictionary
order preserving isomorphisms between finitely generated essential
right ideals of $A^*$.

The monoids ${\sf riAut}(k)$ and ${\sf riAut}_{\sf dict}(k)$
are not groups but they nevertheless have interesting, sometimes 
surprising properties. 
In summary, ${\sf riAut}(k)$ maps homomorphically onto $G_{k,1}$,
and $G_{k,1}$ is the maximum group-homomorphic image. Similarly, 
${\sf riAut}_{\sf dict}(k)$ maps onto $F_{k,1}$, and the latter is
the maximum group-homomorphic image.
Both ${\sf riAut}(k)$ and ${\sf riAut}_{\sf dict}(k)$ are finitely
generated as monoids, their word problem is in {\sf P}, they are 
residually finite, and they are F-inverse. 

Another way to obtain $G_{k,1}$ as a homomorphic image of a residually 
finite F-inverse monoid is to take the $\Gamma$-generated {\it suffix 
expansion} $(G_{k,1}^{\sim {\cal L}})_{_{\Gamma}}$ of $G_{k,1}$, where 
$\Gamma$ is a generating set of $G_{k,1}$.  This will be defined in 
Section 4; the suffix and the prefix expansions were introduced in the 
early 1980's and had a priori no special connection with the 
Thompson-Higman groups. 
We show that $(G_{k,1}^{\sim {\cal L}})_{_{\Gamma}}$ maps 
homomorphically into ${\sf riAut}(k)$, that the map is finite-to-one,
and that it is surjective for certain (finite and infinite) choices of 
$\Gamma$. Thus, the relationship between ${\sf riAut}(k)$ and the prefix
expansion $(G_{k,1}^{\sim {\cal L}})_{_{\Gamma}}$ reveals that certain 
finite generating sets $\Gamma$ of $G_{k,1}$ have special properties. 
In combinatorial group theory it is very rare that some finite 
generating sets of a group behave very differently than others (in 
non-trivial ways). 
Similar results hold for $F_{k,1}$.                 

\bigskip

\noindent {\bf Background definitions and facts} 

\medskip 

We will define the Thompson-Higman groups $G_{k,1}$ and $F_{k,1}$,
as well as the inverse monoids ${\sf riAut}(k)$ and
${\sf riAut}_{\sf dict}(k)$ but we need some preliminary definitions;
we follow \cite{BiThomps} and \cite{BiThomMon} (and indirectly 
\cite{ESc}).

Let $A = \{a_1, \dots, a_k\}$ be a finite alphabet of cardinality
$|A| = k$. The free monoid over $A$ (consisting of all finite strings
over $A$) is denoted by $A^*$. The length of $w \in A^*$ is denoted by
$|w|$, and the empty word is denoted by $\varepsilon$, where
$|\varepsilon| = 0$. The concatenation of $x, y \in A^*$ is denoted by
$xy$ or by $x \cdot y$, and for $B, C \subseteq A^*$ the concatenation
is defined by $BC = \{ xy : x \in B, y \in C\}$.
For $x, y \in A^*$ we say that $x$ is a {\it prefix} of $y$ iff $xs = y$
for some $s \in A^*$. For $S \subseteq A^*$, ${\sf pref}(S)$ denotes
the set of all prefixes of the strings in $S$, including the elements of 
$S$ and $\varepsilon$. We say that 
$x, y \in A^*$ are {\it prefix-comparable} iff $x \in {\sf pref}(y)$ or
$y \in {\sf pref}(x)$. 
A {\it prefix code} is a subset $C \subseteq A^*$ whose elements are
two-by-two prefix-incomparable.
A prefix code $C$ is maximal iff $C$ is not a strict subset of any 
other prefix code.

A set $R \subseteq A^*$ is called a {\it right ideal} iff
$RA^* \subseteq R$, and $R$ is called an {\it essential} right ideal iff 
$R$ intersects every right ideal of $A^*$. 
We say that a right ideal $R$ is {\it generated} by a set 
$C \subseteq A^*$ iff $R$ is the intersection of all right ideals that 
contain $C$; equivalently, $R = CA^*$.
One can prove that a right ideal $R$ has a unique minimal (under 
inclusion) generating set, and that this minimal generating set is a
{\it prefix code}, and that this prefix code is maximal iff $R$ is an 
essential right ideal.
Here we will only consider {\it finitely generated} right ideals, and 
finite prefix codes. 

For a partial function $f:A^* \to A^*$ the domain is denoted by
${\sf Dom}(f)$ and the image by ${\sf Im}(f)$.
A {\em right ideal homomorphism} of $A^*$ is a function 
$\varphi: R_1 \to A^*$ such that ${\sf Dom}(\varphi) = R_1$ is a right 
ideal, and for all $x_1 \in R_1$ and all $w \in A^*$: 
 \ $\varphi(x_1w) = \varphi(x_1) \ w$.
Then one can prove that ${\sf Im}(\varphi)$ is also a right ideal, and 
if $R_1$ is a finitely generated right ideal then ${\sf Im}(\varphi)$ is 
also finitely generated.
We write the action of partial functions on the left of the argument;
equivalently, functions are composed from right to left.
When $\varphi: {\sf Dom}(f) \to A^*$ is injective we call $\varphi$ a 
{\em right ideal isomorphism} between the right ideals ${\sf Dom}(f)$ 
and ${\sf Im}(\varphi)$. 

In this paper we only deal with right ideal isomorphisms for which both
${\sf Dom}(\varphi)$ and ${\sf Im}(\varphi)$ are essential, i.e., their 
prefix codes are maximal. We call such an isomorphism a {\em right ideal 
automorphism of} $A^*$. This does {\em not} mean that 
${\sf Dom}(\varphi) = {\sf Im}(\varphi)$; however, ${\sf Dom}(\varphi)$ 
and ${\sf Im}(\varphi)$ are ``essentially equal'', in the sense that 
every ideal that intersects one also intersects the other and vice versa;
equivalently, ${\sf Dom}(\varphi)$ and ${\sf Im}(\varphi)$ have the same
{\it ends} (see section 1 of \cite{BiRL}). 

A right ideal automorphism $\varphi: R_1 \to R_2$ is uniquely 
determined by its restriction $P_1 \to P_2$, where $P_i$ is the finite 
maximal prefix code that generates $R_i$ ($i = 1,2$). This finite 
bijection $P_1 \to P_2$ is called the {\it table of} $\varphi$. The 
prefix code $P_1$ is called the {\em domain code} of $\varphi$ and is 
denoted by ${\sf domC}(\varphi)$; $P_2$ is called the {\em image code},
denoted by ${\sf imC}(\varphi)$. We can write the table of $\varphi$ 
in the form $\{(x_1,y_1), \ldots, (x_n, y_n)\}$
where $\{x_1, \ldots, x_n\} = P_1$ and
$\{y_1, \ldots, y_n\} = P_2$.

The set of right ideal automorphisms with finite tables, called
${\sf riAut}(k)$, is closed under composition, and the identity map 
on $A^*$ serves as an identity for this multiplication; hence, 
${\sf riAut}(k)$ is a monoid.
Every $\varphi \in {\sf riAut}(k)$ is injective, and its
inverse $\varphi^{-1}: {\sf Im}(\varphi) \to {\sf Dom}(\varphi)$
belongs to ${\sf riAut}(k)$, and satisfies 
 \ $\varphi \circ \varphi^{-1} \circ \varphi = \varphi$, and 
 \ $\varphi^{-1} \circ \varphi \circ \varphi^{-1} = \varphi^{-1}$.
Moreover, $\varphi^{-1}$ is the only element $\psi \in {\sf riAut}(k)$
that satisfies $\varphi \circ \psi \circ \varphi = \varphi$ and
$\psi \circ \varphi \circ \psi = \psi$, by injectiveness of $\varphi$
and $\psi$. Hence, ${\sf riAut}(k)$ is an {\it inverse monoid}. We have
$\varphi^{-1} \circ \varphi = {\sf id}_{{\sf Dom}(\varphi)}$ 
(i.e., the identity map, restricted to ${\sf Dom}(\varphi)$),  and
$\varphi \circ \varphi^{-1}  = {\sf id}_{{\sf Im}(\varphi)}$.
Recall again that in this paper, maps act on the left.

An interesting submonoid of ${\sf riAut}(k)$ is 
${\sf riAut}_{\sf dict}(k)$, the {\it dictionary order preserving} 
automorphisms.  The dictionary order $\leq_{\sf dict}$ is a 
well-order on $A^*$ derived from an order $a_1 < \ldots < a_k$ of 
the alphabet $A$: For $u, v \in A^*$ the dictionary order between them 
is the prefix order, if $u$ and $v$ are prefix-comparable. 
If $u$, $v$ are not prefix-comparable we can write $u = pa_ix$ 
and $v = p a_j y$ where $p \in A^*$ is the longest common prefix of $u$ 
and $v$, where $a_i, a_j \in A$ with $a_i \neq a_j$, and $x, y \in A^*$; 
then the dictionary order between $u$ and $v$ is the same as the order 
between $a_i$ and $a_j$.
An injective partial function $f: A^* \to A^*$ is 
{\it dictionary order preserving} iff for all $u, v \in {\sf Dom}(f):$
  \ $u \leq_{\sf dict} v$ implies $f(u) \leq_{\sf dict} f(v)$.  
One easily proves that if $f$ is dictionary order preserving then
$f^{-1}$ is also dictionary order preserving. The composition of 
dictionary order preserving maps yields a dictionary order preserving
map. Thus ${\sf riAut}_{\sf dict}(k)$ is an inverse submonoid
of ${\sf riAut}(k)$.

\smallskip

We now proceed to the definition of the Thompson-Higman groups 
$G_{k,1}$ and $F_{k,1}$. For a right ideal automorphism 
$\varphi: R_1 \to R_2$, an {\it essential restriction} of $\varphi$
is a right ideal automorphism $\Phi: R'_1 \to R'_2$ such that 
$R'_1, R'_2$ are finitely generated right ideals with 
$R'_1 \subseteq R_1$ and $R'_2 \subseteq R_2$. We also say that 
$\varphi$ is an {\it essential extension} of $\Phi$.
Thompson \cite{Th} (see also \cite{ESc} and \cite{BiThomps}) proved that 
every $\varphi \in {\sf riAut}(k)$ has a {\it unique} maximal essential 
extension in ${\sf riAut}(k)$; we call it ${\sf max}(\varphi)$. 
He showed that an essential restriction
(and, inversely, an essential extension) can be obtained by a finite
number of steps of the following form: In the table 
$\{(x_1,y_1), \ldots, (x_n, y_n)\}$ of $\varphi$, replace some entry 
$(x_i, y_i)$ by  \ $\{ (x_i a_j, y_i a_j) : a_j \in A\}$. 
It follows that ${\sf max}(\varphi) \in {\sf riAut}_{\sf dict}(k)$ 
if $\varphi \in {\sf riAut}_{\sf dict}(k)$. 
We now define $G_{k,1}$ by 

\medskip

\hspace{1.2in}
$G_{k,1} \ = \ \{{\sf max}(\varphi) : \varphi \in {\sf riAut}(k)\}$

\smallskip

\noindent with multiplication

\hspace{1.2in} $\psi \cdot \varphi \ = \ {\sf max}( \psi \circ \varphi)$ .

\smallskip

\noindent This multiplication is associative and turns $G_{k,1}$ into
a group whose identity is the identity map on $A^*$. 
The map \ $\eta: \varphi \mapsto {\sf max}(\varphi)$ 
 \ is a homomorphism from ${\sf riAut}(k)$ onto $G_{k,1}$ (see 
\cite{BiThomps}).  
We define $F_{k,1}$ by

\smallskip

\hspace{1.2in} $F_{k,1} \ = \ $
$\{{\sf max}(\varphi) : \varphi \in {\sf riAut}_{\sf dict}(k)\}$ . 

\smallskip

\noindent Then $F_{k,1}$ is a subgroup of $G_{k,1}$ and it is the 
homomorphic image of ${\sf riAut}_{\sf dict}(k)$ by $\eta$.


\section{Basic properties of ${\sf riAut}(k)$ and
${\sf riAut}_{\sf dict}(k)$ }

We need some more facts about prefix codes.  
Recall that for any set $S \subseteq A^*$, ${\sf pref}(S)$ is the set 
of all prefixes of strings in $S$.
When $S$ is a prefix code we also define 
 \ ${\sf spref}(S) = {\sf pref}(S) \smallsetminus  S$ , 
i.e., the set of {\it strict} prefixes of the strings in $S$.

With any prefix code $P \subset A^*$ one can associate a rooted tree,
called the {\it prefix tree} of $P$, whose vertex set is ${\sf pref}(P)$,
whose edge set is  $\{ (v, va_i) : a_i \in A, \  va_i \in {\sf pref}(P)\}$,
and whose root is $\varepsilon$. The set of leaves of this tree is
$P$ itself.  The non-leaves are called {\it inner vertices} of the
prefix tree, so this set is ${\sf spref}(P)$.
The subtree spanned by ${\sf spref}(P)$ is called the {\it inner tree} of
$P$, and the leaves of the inner tree are called the {\it inner leaves}.

It is well-known that for any {\it maximal} prefix code $P \subset A^*$, 
 \ $|P| = 1 + (k-1) \cdot i$, where $i$
is the number of inner vertices of the prefix tree of $P$ (see
e.g.\ Lemma 6.1(0) in \cite{BicoNP} for a proof and references).
Conversely, for any integer $i \geq 0$ there exists a maximal
prefix code $P \subset A^*$ such that $|P| = 1 + (k-1) \cdot i$.
To summarize:  

\begin{lem} \label{maxPrefCodes} \    
Let $A$ be a finite alphabet with $|A| = k$, let $P \subseteq A^*$ be 
any finite {\em maximal} prefix code, let $V = {\sf pref}(P)$ be the set 
of vertices of the prefix tree of $P$, and let $i$ be the number of inner 
vertices the prefix tree. Then the elements of $P$ are the leaves of the
prefix tree, so \ $|V| = i + |P|$. Moreover,

\medskip

 \ \ \ $|P| \ = \ 1 + (k-1) \cdot i$  ,  

\medskip

 \ \ \ $|V| \ = \ 1 + ki \ = \ \frac{k}{k-1} \cdot |P|$ .  
 \hspace{0.6in} $\Box$ 
\end{lem}

\begin{lem} \label{propertiesRightIdeals1} \ 
If $P_1, P_2 \subset A^*$ are finite maximal prefix codes and if
$P_2 A^* \subseteq P_1 A^*$ then $|P_2| \geq |P_1|$. 
\end{lem}
{\bf Proof.} If $P_2 A^* \subseteq P_1 A^*$, the prefix tree of $P_1$ is 
contained in the prefix tree of $P_2$, hence the prefix tree of $P_2$ has
at least as many vertices as as the prefix tree of $P_1$. 
For maximal prefix codes, the number of leaves grows monotonically with 
the number of vertices (by Lemma \ref{maxPrefCodes}), so the prefix tree 
of $P_2$ has at least as many leaves as the prefix tree of $P_1$.   
  \ \ \ $\Box$

\begin{lem} \label{propertiesRightIdeals2} \  
For all $\varphi, \psi \in {\sf riAut}(k):$
 \ \  $|{\sf domC}(\psi \circ \varphi)| \ \geq \ $ 
${\sf max}\{ |{\sf domC}(\varphi)|, \ |{\sf domC}(\psi)| \}$. 
\end{lem}
{\bf Proof.} For any $\varphi, \psi \in {\sf riAut}(k)$ we have
${\sf Dom}(\psi \varphi) $   $ \subseteq $
${\sf Dom}(\varphi) \cap \varphi^{-1}({\sf Dom}(\psi))$.
Since ${\sf Dom}(\psi \varphi) \subseteq {\sf Dom}(\varphi)$ it follows
(by Lemma \ref{propertiesRightIdeals1})
that $|{\sf domC}(\psi \varphi)| \geq |{\sf domC}(\varphi)|$.
And since ${\sf Dom}(\psi \varphi) \subseteq \varphi^{-1}({\sf Dom}(\psi))$
it follows that
$|{\sf domC}(\psi \varphi)| \geq |\varphi^{-1}({\sf domC}(\psi))|$
$ = |{\sf domC}(\psi)|$. 
 \ \ \ $\Box$ 

\medskip

In semigroup theory the {\it Green relations} play an important role
in the structure of a semigroup. In a monoid $M$, the Green relations
$\leq_{\cal J}$, $\leq_{\cal R}$, and $\leq_{\cal L}$ are preorders,
defined as follows (see \cite{CliffPres, Grillet, Lawson} for more
details). For $s, t \in M$, we have $t \leq_{\cal J} s$ iff
$t = xsy$ for some $x, y \in M$; equivalently, every ideal that
contains $s$ also contains $t$. Similarly, we have $t \leq_{\cal R} s$
iff $t = sy$ for some $y \in M$, and $t \leq_{\cal L} s$ iff $t = xs$
for some $x \in M$.

\begin{pro} \label{Jorder} \ 
The ${\cal J}$-orders of ${\sf riAut}(k)$ and 
${\sf riAut}_{\sf dict}(k)$ are as follows, for all 
$\varphi_1, \varphi_2:$  

\smallskip

\hspace{1in}   $\varphi_2 \leq_{\cal J} \varphi_1$ \ \ iff 
 \ \  $|{\sf domC}(\varphi_2)| \ \geq \  |{\sf domC}(\varphi_1)|$ . 
\end{pro}
{\bf Proof.} $[ \Rightarrow ]$ \ If $\varphi_2 \leq_{\cal J} \varphi_1$ 
then $\varphi_2 = \beta \varphi_1 \alpha$ for some 
$\beta, \alpha \in {\sf riAut}(k)$. This implies (by Lemma
\ref{propertiesRightIdeals2}) that 
$|{\sf domC}(\varphi_2)| \geq  |{\sf domC}(\varphi_1)|$.

\smallskip

\noindent $[ \Leftarrow ]$ \ Let $\varphi_i$ have table $P_i \to Q_i$
($i = 1,2$), with $|Q_2| = |P_2| \geq |P_1| = |Q_1|$. 
Since $P_1, Q_1, P_2, Q_2$ are maximal prefix codes, all have 
cardinalities of the form $|P_i| = |Q_i| = 1 + (k-1) n_i$, where 
$k = |A|$ and $n_i$ is the number of inner vertices of $P_i$ (which is
the same as the number of inner vertices of $Q_i$). 

For $\varphi \in {\sf riAut}(k)$ with table $P \to Q$, consider a 
restriction step; this consists of replacing some entry $(x,y)$ in the 
table of $\varphi$ by the set of entries
$\{ (x a_1, y a_1), \ldots, (x a_k, y a_k) \}$.  This is equivalent to
replacing $\varphi$ by ${\sf id}_{Q'} \circ \varphi \circ {\sf id}_{P'}$
where $P' = (P - \{x\}) \cup xA$ and $Q' = (Q  - \{yx\}) \cup yA$.

By applying restriction steps to $\varphi_1$ we obtain 
 \ $\Phi_1 = {\sf id}_{Q'_1} \circ \varphi_1 \circ {\sf id}_{P'_1}$ 
 \ such that $|P'_1| = |Q'_1| = |P_2| = |Q_2|$. 
Let $[P_2 \to P'_1]$ and $[P'_1 \to P_2]$ be the dictionary-order 
preserving bijections from $P_2$ to $P'_1$, respectively  from 
$P'_1$ to $P_2$.   Then \ $\varphi_2$ $ = $
$ \varphi_2 \circ {\sf id}_{P_2A^*}$ $ = $ 
$\varphi_2 \circ [P'_1 \to P_2] \circ \Phi_1^{-1} \circ \Phi_1$
$\circ [P_2 \to P'_1]$. Hence $\varphi_2 \leq_{\cal J} \Phi_1$
$\leq_{\cal J} \varphi_1$.

When $\varphi_1, \varphi_2 \in {\sf riAut}_{\sf dict}(k)$, the 
multipliers ${\sf id}_{P'_1}$, ${\sf id}_{Q'}$, $[P_2 \to P'_1]$,
$[P'_1 \to P_2]$, and $\varphi_2$, used in the proof, belong to
${\sf riAut}_{\sf dict}(k)$, so the result holds for 
${\sf riAut}_{\sf dict}(k)$ too.
  \ \ \ $\Box$

\bigskip

\noindent 
By definition, a semigroup $S$ is {\em finite}-$\cal J$-{\em above} 
iff for each $s \in S$ the set $\{ x \in S : x \geq_{\cal J} s\}$ is 
finite.
In $A^*$ there are only finitely many maximal prefix codes $P$ with a 
given cardinality $n = |P|$; precisely, it is the number of trees of 
degree $\leq k$ with $i = \frac{n - 1}{k-1}$ vertices. 
Thus we have:

\begin{cor} \label{riAutFinJabove} \ 
The monoids ${\sf riAut}(k)$ and ${\sf riAut}_{\sf dict}(k)$ are 
finite-$\cal J$-above.   \ \ \ $\Box$ 
\end{cor}

\begin{cor} \label{riAutinvlimit} \ 
The monoid ${\sf riAut}(k)$ is a projective limit of finite inverse
monoids, hence ${\sf riAut}(k)$  is residually finite.
\end{cor}
{\bf Proof.} \ Every semigroup $S$ is the projective limit of the 
Rees quotients
 \ $R_s = \{x \in S: x \geq_{\cal J} s\} \cup \{0\}$, as $s$ ranges over 
$S$. By definition, the semigroup $R_s$ is the set
$\{x \in S: x \geq_{\cal J} s\}$ with a zero $0$ added. 
The multiplication in $R_s$ is  $x \cdot y = xy$ (product in
$S$) if $xy \geq_{\cal J} s$, and $x \cdot y = 0$ if 
$xy \not\geq_{\cal J} s$. So, $S$ maps homomorphically onto $R_s$ by
mapping the ideal $\{ x \in S : x \not\geq_{\cal J} s\}$ to $0$.
(see \cite{CliffPres,Grillet}).
Clearly, $S$ is finite-$\cal J$-above iff each Rees quotient semigroup
$R_s$ is finite.
  \ \ \ \ \ $\Box$

\medskip

The formula \ $|P| \ = \ 1 + i \cdot (|A| -1)$, and the
characterization of the $\cal J$-order yield the following. 
\begin{cor} \label{JOmegaChain} \ The $\cal J$-classes of ${\sf riAut}(k)$
form an $\omega$-chain
$( J_i: i \in \omega)$, where
$J_i = \{ \varphi : |{\sf domC}(\varphi)| = 1 + i \cdot (|A| -1) \}$
(and similarly for ${\sf riAut}_{\sf dict}(k)$).

All maximal subgroups of the $\cal J$-class $J_i \subset {\sf riAut}(k)$ 
are isomorphic to the symmetric group ${\mathfrak S}_n$, where 
$n = 1 + i \cdot (|A| -1)$.
The group of units of ${\sf riAut}(k)$ is $J_0$, and $|J_0| = 1$.

In ${\sf riAut}_{\sf dict}(k)$ all subgroups are trivial.
 \ \ \ \ \ $\Box$
\end{cor}

The $\cal R$-order of ${\sf riAut}(k)$ corresponds to the inclusion 
relation between finitely generated right ideal. We have for all
$\varphi_1, \varphi_2 \in {\sf riAut}(k):$
 \ $\varphi_2 \leq_{\cal R} \varphi_1$ \ iff
 \ ${\sf Im}(\varphi_2) \subseteq {\sf Im}(\varphi_1)$ \ iff 
 \ ${\sf imC}(\varphi_1) \subseteq {\sf pref}({\sf imC}(\varphi_2))$.

The $\cal L$-order of ${\sf riAut}(k)$ corresponds to the refinement of 
right-congruences on $A^*$; for injective functions, this is equivalent 
to an inclusion of domains, i.e., we have for all
$\varphi_1, \varphi_2 \in {\sf riAut}(k):$
 \ $\varphi_2 \leq_{\cal L} \varphi_1$ \ iff
 \ ${\sf Dom}(\varphi_2) \subseteq {\sf Dom}(\varphi_1)$ \ iff
 \ ${\sf domC}(\varphi_1) \subseteq {\sf pref}({\sf domC}(\varphi_2))$. 

The set of {\em idempotents} of ${\sf riAut}(k)$ is the set of partial
identities ${\sf id}_{PA^*}$ where $P$ ranges over all maximal finite
prefix codes.
Hence, $\eta^{-1}({\bf 1})$ is the set of idempotents of 
${\sf riAut}(k)$.

\medskip

For a semigroup $S$ a {\em group homomorphism} is, by definition, any 
homomorphism from $S$ onto a group. A group homomorphism 
$h_0: S \twoheadrightarrow G_0$ is called {\it maximum} iff every 
group homomorphism $h: S \twoheadrightarrow G$ factors through $h_0$.
For every inverse semigroup $S$, a maximum group homomorphism $h_0$
exists; $h_0$ is unique, except that $G_0$ is only determined up to
isomorphism. The congruence on $S$ determined by $h_0$ is unique
(see \cite{Grillet,CliffPres}). 

\smallskip

In an inverse semigroup $S$ the {\it natural partial oder} is defined 
by $t \leq s$ iff there exist idempotents $e, e'$ such that 
$t = se = e's$.   
A semigroup $S$ is called F-{\it inverse} iff  $S$ is inverse, 
and every congruence class of the maximum group homomorphism of $S$ has 
exactly one maximum element (in the natural order).  
The uniqueness of maximum essential extension of right ideal 
isomorphisms of $A^*$ means that 
${\sf riAut}(k)$ and ${\sf riAut}_{\sf dict}(k)$ are {\rm F}-inverse.

\begin {pro} \label{maxGrImFinv} \ 
The map \ $\eta: {\sf riAut}(k) \twoheadrightarrow G_{k,1}$ and its
restriction ${\sf riAut}_{\sf dict}(k) \twoheadrightarrow F_{k,1}$ 
are maximum group homomorphisms.
\end{pro}
{\bf Proof.} (1) \ Let $h: {\sf riAut}(k) \twoheadrightarrow G$ be any
group homomorphism. We want to show that for any 
$h: {\sf riAut}(k) \twoheadrightarrow H$, if 
$\eta(\varphi) = \eta(\psi)$ then $h(\varphi) = h(\psi)$.
Let $P = {\sf domC}(\varphi)$, $P' = {\sf imC}(\varphi)$, 
$Q = {\sf domC}(\psi)$, $Q' = {\sf imC}(\psi)$.
The assumption $\eta(\varphi) = \eta(\psi)$ is equivalent to 
${\sf max}(\varphi) = {\sf max}(\psi)$; let $g$ be this element of
$G_{k,1}$.  Then 
 \ $\varphi = {\sf id}_{P'A^*} \circ g \circ {\sf id}_{PA^*}$ \ and
 \ $\psi  = {\sf id}_{Q'A^*} \circ g \circ {\sf id}_{QA^*}$. 
Hence $h(\varphi) = h(g) = h(\psi)$, since $h$ maps idempotents to
the identity of $H$. 
 \ \ \ $\Box$

\bigskip

Besides the maximum group homomorphism 
$\eta: {\sf riAut}(k) \twoheadrightarrow G_{k,1}$ 
there are other homomorphisms on ${\sf riAut}(k)$.
For example, for any $i > 0$ let us define $\eta_i$ as the identity map 
on the $\cal J$-classes $J_0, \ \ldots, \ J_{i-1}$ of ${\sf riAut}(k)$, 
and let $\eta_i$ be defined to be $\eta$ on all $J_j$ for $j \geq i$. 
Then the image monoid of $\eta_i$ is
 \ ${\sf riAut}(k)_i \ = \ $
 \ $J_0 \ \cup \ \ldots \ \cup J_{i-1} \ \cup \ G'_{k,1}$,  where $G'_{k,1}$
is an isomorphic copy of $G_{k,1}$, disjoint from 
$J_0 \cup \ldots \cup J_{i-1}$.
This leads to an $\omega$-chain of homomorphisms 
 
\medskip

 \ \ \ \ \ \ ${\sf riAut}(k) \ \twoheadrightarrow $
 \ ${\sf riAut}(k)_0 \ \twoheadrightarrow$ 
 \ ${\sf riAut}(k)_1 \ \twoheadrightarrow$
 \ $ \ \ldots \ \        \twoheadrightarrow$
 \ ${\sf riAut}(k)_i \ \twoheadrightarrow$
 \ ${\sf riAut}(k)_{i+1} \ \twoheadrightarrow$
 \ $ \ \ldots \ \ $ ; 

\medskip

\noindent the direct limit of this chain is $G_{k,1}$.  
Correspondingly, there exists an $\omega$-chain of progressively finer
congruences on ${\sf riAut}(k)$ whose union is the congruence of the
maximum group homomorphism.


\section{Finite generation }

Higman's method for proving finite generation of $G_{k,1}$ 
(\cite{Hig74} p.\ 24-28) can be adapted to prove the following.

\begin{thm} \label{finGen} \ 
The monoids ${\sf riAut}(k)$ and ${\sf riAut}_{\sf dict}(k)$ are 
finitely generated.
\end{thm}

A remark before we prove the Theorem:
We saw in Lemma \ref{maxPrefCodes} that for any maximal prefix code 
$P \subset A^*$ we have \ $|P| = 1 + (k-1) \cdot i$, 
where $i$ is the number of inner vertices of the prefix tree of $P$.
Conversely, for any $i \geq $ there exists a maximal prefix code 
$P \subset A^*$ such that $|P| = 1 + (k-1) \cdot i$. 
It follows that for all $i \geq 3$ there exists a maximal prefix code 
with at least two inner {\em leaves}. This means that for any $n$ of 
the form $n = 1 + (k-1) \cdot i$ with $i \geq 3$ there exists a prefix
code $P$ of the form \ $P = \{ r a_1,  \ldots, r a_k,$
$ s a_1, \ldots, s a_k, z_{2k+1}, \ldots, z_n\}$. 
E.g., the tree with set of inner vertices 
$\{\varepsilon, \, a_1, \, a_2\}$ has inner leaves 
$a_1$ and $a_2$; the corresponding maximal prefix code is 
 \ $\{a_1, \, a_2\} \, A \ \cup \ \{a_3, \ldots, a_k\}$. 

\bigskip

\noindent {\bf Proof that ${\sf riAut}(k)$ is finitely generated:} 
 \ The following Lemma provides a finite generating set.

\begin{lem} \label{genByInnerThreeriAut} \  
The monoid ${\sf riAut}(k)$ is generated by the set of elements of
${\sf riAut}(k)$ whose domain codes have prefix trees with $\leq 3$ 
inner vertices.
\end{lem}
{\bf Proof.} \ Let $\varphi$ be an element of ${\sf riAut}(k)$ with
table $\{(x_1,y_1), \ldots, (x_n,y_n)\}$, where
$n = 1 + (k-1) \cdot i$ with $i \geq 4$.
Since $i \geq 4 >0$ the prefix codes 
${\sf domC}(\varphi) = \{x_1, \ldots, x_n\}$ and 
${\sf imC}(\varphi) = \{y_1, \ldots, y_n\}$ each have at least one inner 
leaf in their respective prefix tree. Hence ${\sf domC}(\varphi)$ is of 
the form $\{ua_1, \ldots, ua_k, \ x_{i_{k+1}}, \ \ldots, \ x_{i_n} \}$ 
$ = \{x_1, \ldots, x_n\}$, and ${\sf imC}(\varphi)$ is of the form 
$\{v a_1, \ldots, v a_k, \ y_{j_{k+1}}, \ \ldots, \ y_{j_n}\}$ 
$ = \{y_1, \ldots, y_n\}$.

We say that {\it the positions of $\{ua_1, \ldots, ua_k\}$ and
$\{v a_1, \ldots, v a_k\}$ in the table of $\varphi$ overlap} iff
the table contains $(ua_i, v a_j)$ for some $i,j$. 

\medskip

\noindent {\sf Case 1:} \ The positions of $\{ua_1, \ldots, ua_k\}$ and 
$\{v a_1, \ldots, v a_k\}$ in the table of $\varphi$ do {\em not}
overlap.

\smallskip

\noindent Then (for some ordering of the columns) the table of $\varphi$ 
has the form 

\bigskip

$\left[ \begin{array}{lll lll lll }
ua_1 & \ldots & ua_k & x_{k+1} & \ldots & x_{2k} & x_{2k+1} & \ldots
                                                                & x_n \\
y_1 &  \ldots & y_k  & va_1 &  \ldots &   va_k &   y_{2k+1} & \ldots
                                                                & y_n
\end{array} \right]. $

\bigskip

\noindent If $i \geq 3$ there exists a maximal prefix code with at least 
two inner leaves, of the form  

\smallskip

 \ \ \ \ \   
$\{ra_1, \ \ldots, \ ra_k, \ sa_1, \ \ldots, \ sa_k, \ z_{2k+1}, \ $
$ \ldots, \ z_n\}$.  

\smallskip

\noindent We can insert this prefix code as a row into the table of 
$\varphi$, yielding

\bigskip

$\left[ \begin{array}{lll lll lll }
ua_1 & \ldots & ua_k & x_{k+1} & \ldots & x_{2k} & x_{2k+1} & \ldots  
                                                                & x_n \\ 
ra_1 & \ldots & ra_k & sa_1 & \ldots &    sa_k &   z_{2k+1} & \ldots 
                                                                & z_n \\   
y_1 &  \ldots & y_k  & va_1 &  \ldots &   va_k &   y_{2k+1} & \ldots           
                                                                & y_n
\end{array} \right]. $

\bigskip

\noindent This three-row table corresponds to a factorization 
$\varphi = \varphi_2 \circ \varphi_1$ where 

\bigskip 

$\varphi_1 \ = \ $
$\left[ \begin{array}{lll lll lll}
ua_1 & \ldots & ua_k & x_{k+1} & \ldots & x_{2k} & x_{2k+1} & \ldots          
                                                                & x_n \\
ra_1 & \ldots & ra_k & sa_1 & \ldots &    sa_k &   z_{2k+1} & \ldots           
                                                                & z_n 
\end{array} \right],$ 

\bigskip
  
$\varphi_2  \ = \ $
$\left[ \begin{array}{lll lll lll }
ra_1 & \ldots & ra_k & sa_1 & \ldots &    sa_k &   z_{2k+1} & \ldots           
                                                                & z_n \\ 
y_1 &  \ldots & y_k  & va_1 &  \ldots &   va_k &   y_{2k+1} & \ldots
                                                                & y_n
\end{array} \right]. $

\bigskip

\noindent Then we also have the factorization 
$\varphi = \psi_2 \circ \psi_1$ where 

\bigskip

$\psi_1  \ = \ $
$\left[ \begin{array}{lll lll lll }
u & x_{k+1} & \ldots & x_{2k} & x_{2k+1} & \ldots & x_n \\
r & sa_1 & \ldots &    sa_k &   z_{2k+1} & \ldots & z_n
 \end{array} \right],$ 

\bigskip

$\psi_2 \ = \ $
$\left[ \begin{array}{lll lll lll}
ra_1 & \ldots & ra_k  & s & z_{2k+1} & \ldots & z_n \\
y_1 &  \ldots & y_k   & v & y_{2k+1} & \ldots & y_n
 \end{array} \right]. $
  
\bigskip

\noindent {\sf Case 2:} \ The positions of $\{ua_1, \ldots, ua_k\}$ and
$\{v a_1, \ldots, v a_k\}$ in the table of $\varphi$  overlap.

\smallskip

\noindent Then the table of $\varphi$ has the form 

\bigskip

$\left[ \begin{array}{lll lll lll lll}
ua_1 & \ldots & . & \ldots & ua_k & \ldots & . & \ldots & . & \ldots &
                                                           . & \ldots \\
. & \ldots & va_{i_1} & \ldots & . & \ldots & va_{i_k} & \ldots & . 
                                                  & \ldots & . & \ldots
\end{array} \right], $

\bigskip

\noindent where $ua_i$ is in column $i$ and row 1, and $v a_{i_j}$ is
in column $i_j$ and row 2. 

If $i \geq 4$ then a little calculation shows that for all $k \geq 2$: 
 \ \ $n \ = \ 1 + (k-1) \, i \ \geq \ 3k-1$. The fact that $n \geq 3k-1$ 
means that there are at least $k$ columns in the table, in addition to 
the $\leq 2k-1$ columns occupied by $\{ua_1, \ldots, ua_k\}$ 
and $\{v a_1, \ldots, v a_k\}$.  
So we can insert two new rows, each corresponding to a prefix code with 
two inner leaves, as follows: 

\bigskip

$\left[ \begin{array}{lll lll lll lll}
ua_1 & \ldots & . & \ldots & ua_k & \ldots & . & \ldots & . & \ldots & 
                                                           . & \ldots \\
ra_1 & \ldots & . & \ldots & ra_k & \ldots & . & \ldots & sa_1 & \ldots &  
                                                         sa_k & \ldots \\ 
. & \ldots & ra_{i_1} & \ldots & . & \ldots & ra_{i_k} & \ldots &  sa_1 & 
                                              \ldots &  sa_k & \ldots \\ 
. & \ldots & va_{i_1} & \ldots & . & \ldots & va_{i_k} & \ldots & . & 
                                                \ldots & . & \ldots 
\end{array} \right]. $

\bigskip

\noindent This four-row table corresponds to a factorization
$\varphi = \varphi_3' \circ \varphi_2' \circ \varphi_1'$ where

\bigskip

$\varphi_1' \ = \ $
$\left[ \begin{array}{lll lll lll lll}
ua_1 & \ldots & ua_k & \ldots & \ldots & . & \ldots & . & \ldots \\
ra_1 & \ldots & ra_k & \ldots & \ldots & sa_1 & \ldots &
                                                         sa_k & \ldots
\end{array} \right],$ 

\bigskip

$\varphi_2' \ = \ $
$\left[ \begin{array}{lll lll lll lll}
ra_1 & \ldots & . & \ldots & ra_k & \ldots & . & \ldots & sa_1 & \ldots 
                                                      & sa_k & \ldots \\
. & \ldots & ra_{i_1} & \ldots & . & \ldots & ra_{i_k} & \ldots &  
                               sa_1 & \ldots &  sa_k & \ldots 
\end{array} \right],$ \ \ and

\bigskip

$\varphi_3' \ = \ $
$\left[ \begin{array}{lll lll lll lll}
. & \ldots & ra_{i_1} & \ldots & ra_{i_k} & \ldots &  
                            sa_1 & \ldots &  sa_k & \ldots \\
. & \ldots & va_{i_1} & \ldots & va_{i_k} & \ldots & . & \ldots & . & \ldots
\end{array} \right]. $

\bigskip

\noindent Then we also have the factorization 
$\varphi = \psi_3' \circ \psi_2' \circ \psi_1'$ where 

\bigskip

$\psi_1'  \ = \ $
$\left[ \begin{array}{lll lll lll lll}
u & \ldots & \ldots & . & \ldots & . & \ldots \\
r & \ldots & \ldots & sa_1 & \ldots & sa_k & \ldots
\end{array} \right],$

\bigskip

$\psi_2'  \ = \ $
$\left[ \begin{array}{llllllllllll}
ra_{i_1} & \ldots & . & \ldots & ra_{i_k} & \ldots & . & \ldots & s & \ldots \\
. & \ldots & ra_{i_1} & \ldots & . & \ldots & ra_{i_k} & \ldots & s & \ldots
\end{array} \right],$ \ \ \ and

\bigskip

$\psi_3' \ = \ $
$\left[ \begin{array}{llllllllllll}
. & \ldots & r & \ldots &  sa_1 & \ldots &  sa_k & \ldots \\
. & \ldots & v & \ldots & . & \ldots & . & \ldots
\end{array} \right]. $

\bigskip

\noindent In both cases 1 and 2 the factors $\psi_1$, $\psi_2$, 
$\psi_1'$, $\psi_2'$ and $\psi_3'$ of $\varphi$ have tables that have
fewer columns than the table of $\varphi$. We conclude, by induction, 
that every element $\varphi \in {\sf riAut}(k)$ can be written as a 
composition of elements of table-size $< 1 + 4(k-1)$. Hence the table-size 
of these elements will be $\leq 1 + 3(k-1)$ since a maximal prefix code 
has a size of the form $1 + i(k-1)$.     
 \ \ \ \ \ $\Box$

\bigskip

There are only finitely many elements in ${\sf riAut}(k)$ with 
table-size $\leq 1 + 3 (k-1)$, so ${\sf riAut}(k)$ is finitely 
generated. This proves Theorem \ref{finGen} for ${\sf riAut}(k)$. 
 \ \ \ \ \ $\Box$

\bigskip

\bigskip


\noindent {\bf Proof that ${\sf riAut}_{\sf dict}(k)$ is finitely 
generated:} 
  \ The following Lemma provides a finite generating set.

\begin{lem} \label{genByInnerThreeriAutdict} \   
The monoid ${\sf riAut}_{\sf dict}(k)$ is generated by the set of 
elements of ${\sf riAut}_{\sf dict}(k)$ whose domain codes have prefix 
trees with $\leq k+1$ inner vertices, i.e., whose table size is
 \ $\leq k^2$ .
\end{lem}
The proof is similar to the proof of Lemma \ref{genByInnerThreeriAut}, 
with the added constraint that all factors must preserve the dictionary 
order. 
To ensure that all elements of ${\sf riAut}(k)$ considered here
preserve the dictionary order, we will write every maximal prefix code as 
a {\it sequence}, according to strictly increasing dictionary order.  
Since the alphabet $A$ is ordered (by $a_1 < \ldots < a_k$) the prefix 
tree of a prefix code is now an {\it oriented tree}, i.e., the set of
children of every vertex is ordered.
An element $\varphi \in {\sf riAut}_{\sf dict}(k)$ has a table 
$\left[ \begin{array}{lll}
 x_1 & \ldots & x_n \\
 y_1 & \ldots & y_n
 \end{array} \right]$ where $x_1 <_{\sf dict} \ldots <_{\sf dict} x_n$
and $y_1 <_{\sf dict} \ldots <_{\sf dict} y_n$. When we insert one or two
rows into the table of $\varphi$, as we did in the proof of Lemma
\ref{genByInnerThreeriAut}, the new rows must also be in increasing 
dictionary order. 

In the proof of Lemma \ref{genByInnerThreeriAut} we used the fact that 
for all $i \geq 3$ there exists a maximal prefix code with at 
least two inner leaves. For ${\sf riAut}_{\sf dict}(k)$ we need some 
control over the {\it position} of these inner leaves (according to the 
dictionary order of leaves of the prefix tree): 

\begin{lem} \label{twoIntLeaves} \ 
Let $P$ be a maximal prefix code, let $z$ be a leaf of the inner 
tree of $P$, let $\ell$ be the number of leaves of $P$ that are 
strictly to the left of $z$, and let $r$ be the number of leaves 
of $P$ that are strictly to the right of $z$.  
In other words, $|P| = \ell + k + r$ and $P$ is of the form 

\smallskip

 \ \ \ \ \ $P \ = \ (x_1, \ \ldots, \ x_{\ell}, \ za_1, \ldots, za_k,$
$ \ x_{\ell + k +1}, \ \ldots, \ x_{\ell + k + r})$.

\smallskip

\noindent Then if \ $|P| > 1 + (k-1) \, (k+1)$ $(= k^2$), there 
exists a maximal prefix code $Q$ such that: 
 
\smallskip

\noindent $\bullet$ \ $|Q| = |P|$ \ ($= \ell + k + r$); 

\smallskip

\noindent $\bullet$ \ $Q$ has an inner leaf $Z$ such that $Q$ has 
$\ell$ leaves strictly to the left of $Z$ and $r$ leaves strictly 
to the right of $Z$; 

\smallskip

\noindent $\bullet$ \ $Q$ has an additional inner leaf $Z'$ 
($\neq Z$). 
\end{lem}
{\bf Proof.}
If $P$ has two inner leaves we can take $Q$ to be $P$ itself. 
Let us assume now that $P$ has only one inner leaf, i.e., the inner
tree of $P$ is just a path; let $z$ be the label of this path.

For reasons that will be clear below (Case 3) we assume that $z$ has 
length $|z| \geq k+1$. This is always the case if the number of 
inner vertices of $P$ is at least $k+2$ (since the inner tree is a 
path), i.e., if $|P| \geq 1 + (k-1) \, (k+2)$. Equivalently (since 
$P$ is a maximal prefix code),  $|P| > 1 + (k-1) \, (k+1)$.

The maximal prefix code $Q$ (with inner leaves $Z$, $Z'$, etc.) is 
constructed from the maximal prefix code $P$ by removing one edge from 
the inner path $z$, reconnecting, and possibly shifting, so as to make 
a new path $Z$ of length $|Z| = |z| - 1$. 
Next, an additional inner leaf is attached at an appropriate place 
on the side of the inner path $Z$. The details are given next.
Note that $|z| \geq k+1$ implies $|z| \geq 3$. 

\medskip

\noindent {\sf Case 1:} \ \ $z$ contains $a_1$ and an additional letter 
$a_j \neq a_k$. 

We have $z = u a_1 v$ for some $u, v \in A^*$.  To construct $Q$ from 
$P$ we remove an edge with label $a_1$ from the path $z$ and reconnect. 
The new path is $Z = uv$; also, $Z = Xa_jY$ with $j<k$, for some 
$X, Y \in A^*$. Since a vertex of the form $ua_1$ has no 
left-siblings, the replacement of $z$ by $Z$ does not change $\ell$; 
but the number of inner vertices has been decreased by 1. 
To preserve $|P|$ we attach an additional child to vertex $X$ on the 
right of $Xa_j$, i.e., we create a new inner vertex $Z' = X a_{j+1}$ 
in $Q$.

\smallskip

\noindent {\sf Case 2:} \ \ $z$ contains $a_k$ and an additional letter 
$a_i \neq a_1$.

This case is left-right symmetric to case 1, since preserving $|P|$ and
$r$ is equivalent to preserving $|P|$ and $\ell$.

\smallskip

\noindent {\sf Case 3:} \ $z$ contains no occurrences of $a_1$ nor $a_k$. 

Then $z$ has the form $a_j a_{i_1} \ldots a_{i_{j-1}} w$ with 
$2 \leq j, \ i_1, \ \ldots , \ i_{j-1} \leq k-1$, and $w \in A^*$; 
 \ recall that $|z| \geq k+1$. 
We remove $a_j$ from $z$. This removes one vertex from
the inner tree (since the inner tree is a path), and decreases $\ell$
by the amount $j-1$.  In order to preserve $\ell$ we let 
 \ $Z = a_{i_1 +1} \ldots a_{i_{j-1} +1} w$.
In order to preserve the total number of inner vertices, and in order to 
create an additional inner leaf $Z'$ we add one inner vertex, namely
$Z' = a_k$.

\smallskip

Let us verify that his completes the construction of $Q$, i.e., that 
cases 1, 2, and 3 exhaust all possibilities. If case 3 does not hold, 
$z$ contains $a_1$ or $a_k$. If $z$ contains $a_1$ but case 1 does not 
hold, $z$ consists of $\geq k+1$ copies of $a_k$; then case 2 holds. 
If $z$ contains $a_k$ but case 2 does not hold, $z$ consists of $\geq k+1$ 
copies of $a_1$; then case 1 holds. 
 \ \ \ \ \ $\Box$

\medskip

\noindent {\bf Proof of Lemma \ref{genByInnerThreeriAutdict}:} \  
Just as in the proof of Lemma \ref{genByInnerThreeriAut},
we can factor any element $\varphi \in {\sf riAut}_{\sf dict}(k)$ 
into a product of elements of ${\sf riAut}_{\sf dict}(k)$ with smaller
tables, whenever the table of $\varphi$ has size 
$> 1 + (k-1) \, (k+1) \ = \ k^2$.
 \ \ \ $\Box$

\bigskip

\noindent {\bf Open problem:} \ Are ${\sf riAut}(k)$ and
${\sf riAut}_{\sf dict}(k)$ finitely presented?


\section{ The suffix expansion } 

For any monoid $M$ let ${\cal P}_{_{\sf FIN}}(M)$ denote the set of 
finite subsets of $M$; the union operation makes 
${\cal P}_{_{\sf FIN}}(M)$ a semilattice. We define the left 
semidirect product $M \ltimes ({\cal P}_{_{\sf FIN}}(M), \cup)$  
with multiplication

\medskip

 \ \ \ \ \ $(y, T) \cdot (x, S) \ = \ (y x , \ T x \cup S)$,

\medskip

\noindent where $T x = \{ t x : t \in T\}$. The subset 
 \ $\{ (x, S) \in M \ltimes {\cal P}_{_{\sf FIN}}(M) \ : \ $
$\{x, 1\} \subseteq S \}$ \ of $M \ltimes {\cal P}_{_{\sf FIN}}(M)$ 
is closed under multiplication.
Let $\Gamma$ be a generating set of $M$.  By definition, the 
{\it $\Gamma$-generated suffix expansion} of $M$, denoted by 
$(\tilde{M}^{\cal L})_{_{\Gamma}}$, is the submonoid of 
$M \ltimes {\cal P}_{_{\sf FIN}}(M)$ generated by the 
subset $\{(\gamma, \{\gamma,1\}) : \gamma \in \Gamma\}$. 
When $\Gamma$ is $M$ itself, \ $(\tilde{M}^{\cal L})_{_{\Gamma}}$ \ is
simply denoted by $\tilde{M}^{\cal L}$.
The monoid $(\tilde{M}^{\cal L})_{_{\Gamma}}$ maps homomorphically onto
$M$ by the projection $(x, S) \mapsto x$. 

The transformation from $M$ to $\tilde{M}^{\cal L}$ can also be applied 
to monoid homomorphisms, and thus becomes a functor. For details, see
\cite{BiRhAlmost} and \cite{BiRhGrvia}, where the suffix expansion was 
introduced and where many properties were proved.
For example, for any group $G$ the inverse monoid 
$(\tilde{G}^{\cal L})_{_{\Gamma}}$ is finite-$\cal J$-above, hence it is
a projective limit of finite inverse monoids \cite{BiRhGrvia}.
The idea of semigroup expansions is due to John Rhodes (see 
\cite{Tilson}).

Dually (switching left and right) one defines the {\it prefix expansion} 
$(\tilde{M}^{\cal R})_{_{\Gamma}}$ as the submonoid generated by 
$\{ (\{1,\gamma\}, \gamma) : \gamma \in \Gamma\}$ in  the right
semidirect product $({\cal P}_{_{\sf FIN}}(M), \cup) \rtimes M$.
For any group $G$, $(\tilde{G}^{\cal L})_{_{\Gamma}}$ and
$(\tilde{G}^{\cal R})_{_{\Gamma}}$ are isomorphic \cite{BiRhGrvia};
in this paper we will only work with the suffix expansion.
Also \cite{Szendrei}, for a group the underlying set of
$\tilde{G}^{\cal L}$ is all of 
 \ $\{ (g, S) : \ g \in G, \ $
$S \in {\cal P}_{_{\sf FIN}}(G), \ \{g,1\} \subseteq S\}$, and similarly
for $\tilde{G}^{\cal R}$.
Szendrei \cite{Szendrei} proved that $(.)^{\sim {\cal R}}$ is a functor 
from the category of groups to the category of F-inverse monoids, and 
that it is the left-adjoint of the maximum-group-image functor (i.e., 
the functor which maps an F-inverse semigroup to its maximum group 
homomorphic image). 

The main result of this section is that for certain generating sets 
$\Gamma$ of $G_{k,1}$, the suffix expansion 
$(G_{k,1}^{\sim {\cal L}})_{_{\Gamma}}$ maps onto ${\sf riAut}(k)$;
similarly, $(F_{k,1}^{\sim {\cal L}})_{_{\Gamma}}$ maps
onto ${\sf riAut}(k)_{\sf dict}$. We need some preliminary results.

\begin{lem} \label{DomComposite} \ 
Let $f_1: X_1 \to X_2$, $f_2: X_2 \to X_3$ be any partial functions.
Then \ \ ${\sf Dom}(f_2 \circ f_1) = $
${\sf Dom}( f_1) \ \cap \ f_1^{-1}({\sf Dom}(f_2))$.
\end{lem}
{\bf Proof.} We have $x \in {\sf Dom}(f_2 \circ f_1)$ \ iff 
 \ $x \in {\sf Dom}(f_1)$ and $f_1(x) \in {\sf Dom}(f_2)$. The latter 
is equivalent to $x \in f_1^{-1}({\sf Dom}(f_2))$.
 \ \ \ $\Box$

\medskip

We first give an embedding of ${\sf riAut}(k)$ into a semidirect 
product of the Thompson-Higman group $G_{k,1}$ and a semilattice.
Here each element of $G_{k,1}$ is represented by a maximally extended
element of ${\sf riAut}(k)$.

Let ${\cal I_R}$ be the set of finitely generated essential right ideals
of $A^*$. Each such ideal is of the form $PA^*$ where $P$ is a finite 
maximal prefix code. One can prove that the intersection of two 
essential right ideals is an essential right ideal (\cite{BiThomps}, 
Lemma A.2, p.\ 608), and that this intersection is finitely generated 
(\cite{BiThomps}, Lemma 3.3, p.\ 579).  
Thus, ${\cal I_R}$ is closed under intersection, so $({\cal I_R}, \cap)$ 
is a semilattice.  We consider the semidirect product 
$G_{k,1} \ltimes ({\cal I_R}, \cap)$ with multiplication 

\smallskip

 \ \ \ \ \ $(g_2, P_2 A^*) * (g_1, P_1 A^*) $
$ \ = \ $ $(g_2 g_1, \ g_1^{-1}(P_2 A^*) \ \cap \ P_1 A^*)$. 

\smallskip

\noindent It is easy to prove that this multiplication is associative. 
This semidirect product projects homomorphically onto $G_{k,1}$.  
Similarly, for the group $F_{k,1}$ we consider the submonoid 
$F_{k,1} \ltimes ({\cal I_R}, \cap)$ of 
$G_{k,1} \ltimes ({\cal I_R}, \cap)$. 

Recall that $\eta$ denotes the surmorphism
${\sf riAut}(k) \twoheadrightarrow G_{k,1}$, or its restriction
${\sf riAut}_{\sf dict}(k) \twoheadrightarrow F_{k,1}$.

\begin{pro} \label{retractriAutofsemi} \  
The monoid ${\sf riAut}(k)$ is a {\em retract} of
$G_{k,1} \ltimes {\cal I_R}$
by the maps 

\smallskip

 \ \ \ \ \ $e: \ \varphi \in {\sf riAut}(k) \ \hookrightarrow $
$ \ (\eta(\varphi), \ {\sf Dom}(\varphi)) $
$ \ \in \ G_{k,1} \ltimes {\cal I_R}$ \ \ (embedding), and 

\smallskip

 \ \ \ \ \ $e': \ (g, PA^*) \in G_{k,1} \ltimes {\cal I_R}$
$ \ \twoheadrightarrow \ g_{_{PA^*}} \in {\sf riAut}(k)$ ,  

\smallskip

\noindent where $g_{_{PA^*}}$ denotes the restriction of $g$ to $PA^*$.  
 \ So $e({\sf riAut}(k))$ is isomorphic to ${\sf riAut}(k)$.

\smallskip

A similar result holds for $F_{k,1}$, namely, 
${\sf riAut}_{\sf dict}(k)$ is a retract of
$F_{k,1} \ltimes {\cal I_R}$, and $e({\sf riAut}_{\sf dict}(k)$
is isomorphic to ${\sf riAut}_{\sf dict}(k)$.  This is obtained by
restricting $e$ to 
${\sf riAut}_{\sf dict}(k) \hookrightarrow F_{k,1} \ltimes {\cal I_R}$ 
and restricting $e'$ to 
$F_{k,1} \ltimes {\cal I_R} \twoheadrightarrow {\sf riAut}_{\sf dict}(k)$.
\end{pro}
{\bf Proof.} \ By the definition of the elements of ${\sf riAut}(k)$ 
the map $e$ is total and injective. That $e$ is a homomorphism follows 
from Lemma \ref{DomComposite} and the fact that 
$\varphi(x) = (\eta(\varphi))(x)$ when $x \in {\sf Dom}(\varphi)$. 
The restriction map 
 \ $e': (g, PA^*) \mapsto g_{_{PA^*}} \in {\sf riAut}(k)$ 
 \ is clearly surjective. It is a homomorphism by Lemma \ref{DomComposite}.
The retraction property obviously holds, namely,  
 \ $(e(\varphi))_{_{{\sf Dom}(\varphi)}} = \varphi$, and $e'$ is injective
on $e({\sf riAut}(k))$. 
 \ \ \ \ \ $\Box$

\medskip

\noindent The map \ $e: {\sf riAut}(k) \to G_{k,1} \ltimes {\cal I_R}$ 
 \ and its restriction 
 \ ${\sf riAut}_{\sf dict}(k) \to F_{k,1} \ltimes {\cal I_R}$
 \ are not surjective; the respective images are 

\smallskip

 \ \ \ $e({\sf riAut}(k))$
$ \ = \ \{ (g, PA^*) \in G_{k,1} \ltimes {\cal I_R}$
$ \ : \ PA^* \subseteq {\sf Dom}(g)\}$ 
 \ \ $ \subsetneqq \ \ G_{k,1} \ltimes {\cal I_R}$ , 

\smallskip

 \ \ \ $e({\sf riAut}_{\sf dict}(k)$
$ \ = \ \{ (g, PA^*) \in F_{k,1} \ltimes {\cal I_R}$
$ \ : \ PA^* \subseteq {\sf Dom}(g)\}$ 
 \ \ $ \subsetneqq \ \ F_{k,1} \ltimes {\cal I_R}$.

\bigskip

We now prove that the $\Gamma$-generated suffix expansions of 
the Thompson-Higman groups $G_{k,1}$ and $F_{k,1}$ map homomorphically
onto ${\sf riAut}(k)$, respectively ${\sf riAut}_{\sf dict}(k)$.
In the case of $G_{k,1}^{\sim {\cal L}}$ and 
$F_{k,1}^{\sim {\cal L}}$ (i.e.\ when $\Gamma = G_{k,1}$ respectively 
$F_{k,1}$), the ``into''-part of Theorem \ref{GtildaontoriAut} follows
from Szendrei's Corollary 3 in \cite{Szendrei}.


\begin{thm} \label{GtildaontoriAut} \   

\noindent {\bf (1)} 
 \ For any generating set $\Gamma$ of the Thompson-Higman group 
$G_{k,1}$ the suffix expansion $(G_{k,1}^{\sim {\cal L}})_{_{\Gamma}}$ 
maps homomorphically into ${\sf riAut}(k)$ by the map

\medskip

 \ \ \ \ \ \     
$\rho: \ \ (g,S) \in (G_{k,1}^{\sim {\cal L}})_{_{\Gamma}}$
$ \ \ \longmapsto \ \ $
$(g, \ \bigcap_{h \in S} {\sf Dom}(h)) \ \in \ e({\sf riAut}(k))$
 \ \ ($\simeq {\sf riAut}(k)$).

\medskip

\noindent Let us assume in addition that $\Gamma$ satisfies the 
following surjectiveness condition: 
There is a generating set $\Delta$ of ${\sf riAut}(k)$ such that 

\smallskip

\hspace{1.1in} 
$(\forall \delta \in \Delta) (\exists \gamma \in \Gamma) :$
 \ ${\sf Dom}(\delta) = {\sf Dom}(\gamma)$. 

\smallskip

\noindent
Then the homomorphism $\rho$ is onto $e({\sf riAut}(k))$, and 
$e' \circ \rho$ is onto ${\sf riAut}(k)$.

\smallskip

\noindent {\bf (2)}
 \ Similarly, if $\Gamma$ is a generating set of the Thompson group 
$F_{k,1}$ then the suffix expansion 
$(F_{k,1}^{\sim {\cal L}})_{_{\Gamma}}$ maps homomorphically into 
${\sf riAut}_{\sf dict}(k)$, by restricting the map $e' \circ \rho$
to $(F_{k,1}^{\sim {\cal L}})_{_{\Gamma}}$.
This map is onto ${\sf riAut}_{\sf dict}(k)$ if $\Gamma$ satisfies 
the condition that there is a generating set $\Delta$ of 
${\sf riAut}_{\sf dict}(k)$ such that 
 \ $(\forall \delta \in \Delta) (\exists \gamma \in \Gamma) :$
 \ ${\sf Dom}(\delta) = {\sf Dom}(\gamma)$.
\end{thm}
{\bf Proof.} {\bf (a)} 
 \ {\it $\rho: (G_{k,1}^{\sim {\cal L}})_{_{\Gamma}}$
$ \ \to \ G_{k,1} \ltimes {\cal I_R}$ \
 is a homomorphism (for any generating set $\Gamma$ of $G_{k,1}$):} 

\smallskip

For 
 \ $(g_2,S_2), \ (g_1,S_1) \in (G_{k,1}^{\sim {\cal L}})_{_{\Gamma}}$
 \ multiplication is defined by 
 \ $(g_2,S_2) \cdot (g_1,S_1) = (g_2 g_1, \ S_2 g_1 \cup S_1)$, thus 

\smallskip

 \ \ \ $\rho((g_2,S_2) \cdot (g_1,S_1)) \ = \  $
$(g_2 g_1, \ \bigcap_{h \in S_2 g_1 \cup S_1} {\sf Dom}(h) \, )$.

\smallskip

\noindent Moreover, 

\smallskip

$\bigcap_{h \in S_2 g_1 \cup S_1} {\sf Dom}(h) \ \ = \ \ $
$\bigcap_{k \in S_2 g_1} {\sf Dom}(k) \ \ \cap \ $
$\bigcap_{h_1 \in S_1} {\sf Dom}(h_1) $

\medskip

$ \ = \ \ \bigcap_{h_2 \in S_2} {\sf Dom}(h_2g_1) \ \ \cap \ $
$\bigcap_{h_1 \in S_1} {\sf Dom}(h_1)$

\medskip

$ \ = \ \ g_1^{-1}(\bigcap_{h_2 \in S_2} {\sf Dom}(h_2)) \ \cap \ $
${\sf Dom}(g_1) \ \cap \ \bigcap_{h_1 \in S_1} {\sf Dom}(h_1)$ ;

\smallskip

\noindent the last equality holds since  
 \ ${\sf Dom}(h_2g_1) = g_1^{-1}({\sf Dom}(h_2)) \ \cap \ {\sf Dom}(g_1)$
 \ (by Lemma \ref{DomComposite}). 
Also, since $g_1 \in S_1$ we have 
 \ ${\sf Dom}(g_1) \ \cap \ \bigcap_{h_1 \in S_1} {\sf Dom}(h_1) \ = \ $
$\bigcap_{h_1 \in S_1} {\sf Dom}(h_1)$.
 \ Thus we have

\smallskip

$\rho((g_2,S_2) \cdot (g_1,S_1)) \ = \ $
$\big( g_2 g_1, \ g_1^{-1}(\bigcap_{h_2 \in S_2} {\sf Dom}(h_2))$
$  \ \cap \ $
$\bigcap_{h_1 \in S_1} {\sf Dom}(h_1) \, \big)$.

\smallskip 

\noindent
A straightforward multiplication in $G_{k,1} \ltimes {\cal I_R}$ shows 
that the latter is also equal to the product
 \ $\rho((g_2,S_2)) * \rho((g_1,S_1))$.

\medskip

\noindent {\bf (b)} \ {\it $\rho$ maps 
$(G_{k,1}^{\sim {\cal L}})_{_{\Gamma}}$ into $e({\sf riAut}(k))$
(for any generating set $\Gamma$ of $G_{k,1}$): } 

\smallskip

We want to show that for every
$(g,S) \in (G_{k,1}^{\sim {\cal L}})_{_{\Gamma}}$, 
 \ $\rho((g,S)) = (g, \bigcap_{h \in S} {\sf Dom}(h))$ \ is equal to
$(g, PA^*)$ for some finite maximal prefix code $P \subset A^*$ such 
that $PA^* \subseteq {\sf Dom}(g)$. 
We saw that the intersection of finitely may finitely generated essential
right ideals is a finitely generated essential right ideal, so 
 \ $\bigcap_{h \in S} {\sf Dom}(h) = PA^*$ \ for some finite maximal 
prefix code $P$.  Moreover, $g \in S$ for every
$(g,S) \in (G_{k,1}^{\sim {\cal L}})_{_{\Gamma}}$, hence ($PA^* =$) 
$\bigcap_{h \in S} {\sf Dom}(h)$ $ \ \subseteq \ {\sf Dom}(g)$.

\medskip

\noindent {\bf (c)} \ {\it $\rho$ maps 
$(G_{k,1}^{\sim {\cal L}})_{_{\Gamma}}$ onto $e({\sf riAut}(k))$ 
(if $\Gamma$ satisfies the surjectiveness condition): }

\smallskip

Let $\Gamma$ be a generating set of $G_{k,1}$ satisfying the condition 
of the Theorem, and let $\Delta$ be a corresponding generating set of 
${\sf riAut}(k)$.
To show that $\rho$ maps onto $e({\sf riAut}(k))$ it is sufficient to
show that $e(\Delta)$ is in the image of $\rho$.

By to the definition of $e$ (in Prop.\ \ref{retractriAutofsemi}), 
for any $\delta \in \Delta$, $e(\delta)$ is of the form 
$e(\delta) = (g, PA^*)$, where $g = \eta(\delta) \in G_{k,1}$ and 
$P = {\sf domC}(\delta)$, with $PA^* \subseteq {\sf Dom}(g)$.
By the condition of the Theorem there exists $\gamma \in \Gamma$ such 
that ${\sf domC}(\gamma) = {\sf domC}(\delta) = P$. 
Then $\gamma^{-1} \circ \gamma = {\sf id}_{PA^*}$, so 
 \ $\delta = g \circ \gamma^{-1} \circ \gamma$ \ (product in 
${\sf riAut}(k)$).  Then, by multiplying in 
$\rho((G_{k,1}^{\sim {\cal L}})_{_{\Gamma}})$ we obtain

\smallskip

$\rho \big( (g,  \{g, {\bf 1}\})  \ \cdot \ $ 
$(\gamma^{-1}, \{\gamma^{-1}, {\bf 1}\}) \ \cdot \ $
$(\gamma , \{\gamma, {\bf 1}\}) \big) $  
$ \ = \ $
$\rho \big( (g,  \{g, {\bf 1}\}) \ \cdot \ $ 
$({\bf 1}, \ \{\gamma^{-1}, {\bf 1}\} \gamma$ $\cup$
                             $ \{\gamma, {\bf 1}\}) \big)$

\smallskip

$ \ = \ $
$\rho \big( (g,  \{g, {\bf 1}\}) \big) \ \cdot \ $
$\rho \big( ({\bf 1}, \{\gamma, {\bf 1}\}) \big)$
$ \ = \ $
$(g, {\sf Dom}(g)) \ \cdot \ ({\bf 1}, PA^*) \ \ = \ \ $
$(g , \ {\bf 1}^{-1}({\sf Dom}(g)) \cap  PA^*)  $

\smallskip

$ \ = \ $
$(g , PA^*)$, \ \ since $PA^* \subseteq {\sf Dom}(g)$.

\smallskip

\noindent Thus, $(g , PA^*) = e(\delta)$ 
$ = \rho \big( (g,  \{g, {\bf 1}\}) \ \cdot \ $
$(\gamma^{-1}, \{\gamma^{-1}, {\bf 1}\}) \ \cdot \ $
$(\gamma , \{\gamma, {\bf 1}\}) \big) \ \in \ $
$\rho((G_{k,1}^{\sim {\cal L}})_{_{\Gamma}})$.  So $e(\delta)$  is in 
$\rho((G_{k,1}^{\sim {\cal L}})_{_{\Gamma}})$ for every 
$\delta \in \Delta$.  

\medskip

\noindent {\bf (d)} \ The same proof applies to 
 \ $(F_{k,1}^{\sim {\cal L}})_{_{\Gamma}} \to $
${\sf riAut}_{\sf dict}(k)$.
  \ \ \ \ \ $\Box$

\medskip

We will prove next that the surjectiveness condition in Theorem
\ref{GtildaontoriAut} holds for some, but not all generating sets 
$\Gamma$, and that it is necessary.  We first need a Lemma. 

\begin{lem} \label{everyPisdomC} \
For every finite maximal prefix code $P \subset A^*$ there is an element
$\varphi \in F_{k,1}$ $(\subset G_{k,1})$ such that 
$P = {\sf domC}(\varphi)$ \ (when $\varphi$ is in maximally extended
form).
\end{lem}
{\bf Proof.} \ Let $|P| = 1 + i \, (k-1)$ where $i$ is the number of
inner vertices of the prefix tree of $P$.  Consider the maximal prefix
code $Q$ whose inner tree consists of the path $a_1^{i-2} a_2$. Hence 
the set of inner vertices of $Q$ is ${\sf pref}(a_1^{i-2} a_2)$, and
$Q$ has only one inner leaf. Also, $Q$ has $i$ inner vertices, so
$|Q| = |P|$.

If $Q \neq P$ then $P$ does not have $a_1^{i-2} a_2$ as an inner leaf.
Indeed, if the inner tree of $P$ is not a path, it will not contain any 
path of length $i-1$; and if $P$ is a path but $P \neq Q$, this path is
different from $a_1^{i-2} a_2$.  
Hence the dictionary-order preserving bijection
$\varphi = (P \to Q) \in F_{k,1}$ is in maximally extended form. Indeed,
extensions steps of an element of $G_{k,1}$ can only happen at a 
common inner leaf of the domain code $P$ and the image code $Q$.
Hence, $P = {\sf domC}(\varphi)$.

If $Q = P$, consider the maximal prefix code $Q'$ whose inner tree
consists of the path $a_1^{i-1}$. Then $P$ does not have $a_1^i$ as an
inner leaf, hence $\varphi' = (P \to Q') \in F_{k,1}$ is in
maximally extended form, so $P = {\sf domC}(\varphi')$.
 \ \ \ \ \ $\Box$

\begin{pro} \label{existenceofGammaG} \ 

\noindent {\bf (1)} 
 \ For every generating set $\Delta$ of ${\sf riAut}(k)$ 
there exists a generating set $\Gamma$ of $G_{k,1}$ that satisfies
the surjectiveness condition of Theorem \ref{GtildaontoriAut} (namely, 
for every $\delta \in \Delta$ there exists $\gamma \in \Gamma$ with 
${\sf Dom}(\delta) = {\sf Dom}(\gamma)$). 

If $\Delta$ is finite then $\Gamma$ is finite and 
 \ $|\Gamma| \leq 2 \cdot |\Delta|$. 

The generating set $\Gamma = G_{k,1}$ also satisfies the surjectiveness
condition.

\smallskip

\noindent {\bf (2)} \ The condition on $\Gamma$ in Theorem
\ref{GtildaontoriAut} is necessary for the surjectiveness of $\rho$,
in general. 

\smallskip

\noindent {\bf (3)} \ Not every generating set $\Gamma$ of $G_{k,1}$
satisfies the surjectiveness condition.
More strongly, 
for some generating set $\Gamma$ of $G_{2,1}$ there is no surjective 
homomorphism from $(G_{k,1}^{\sim {\cal L}})_{_{\Gamma}}$ onto
${\sf riAut}(k)$.
\end{pro}
{\bf Proof.} \ (1) Let $\Delta$ be any generating set of 
${\sf riAut}(k)$.  By Lemma \ref{everyPisdomC}, for each 
$\delta \in \Delta$ there exists $\varphi_{\delta} \in F_{k,1}$ with
${\sf domC}(\varphi_{\delta})  = {\sf domC}(\delta)$.
Let \ $\Gamma \ = \ \eta(\Delta) \cup \{\varphi_{\delta} $
$: \delta \in \Delta\}$. Then $\Gamma$ has the claimed properties.
When $\Delta$ is finite we have $|\Gamma| \leq  $
$|\eta(\Delta)| + |\{\varphi_{\delta} : \delta \in \Delta\}|$  
$ \leq 2 \cdot |\Delta|$. In Section 3 we proved that ${\sf riAut}(k)$
has a finite generating set.

For $\Gamma = G_{k,1}$ and $\Delta = {\sf riAut}(k)$, every finite 
maximal prefix code $P$ occurs as a domain code of an element of 
$G_{k,1}$ and as the domain code of an element of ${\sf riAut}(k)$; 
so the surjectiveness condition of Theorem \ref{GtildaontoriAut} 
applies to $\Gamma$.

\smallskip

\noindent (2) \ For any finite generating set $\Gamma$ of $G_{k,1}$ 
the corresponding generating set of 
$(G_{k,1}^{\sim {\cal L}})_{_{\Gamma}}$ is \ $\tilde{\Gamma} = $
$ \{ (\gamma, \{\gamma,{\bf 1}\}): \gamma \in \Gamma\}$. 
If Theorem \ref{GtildaontoriAut} holds for $\Gamma$, i.e., $\Gamma$ is 
such that 
$\rho: (G_{k,1}^{\sim {\cal L}})_{_{\Gamma}} \to e({\sf riAut}(k))$ 
is surjective, then $\rho(\tilde{\Gamma}) = $
$\{ (\gamma, {\sf Dom}(\gamma)) : \gamma \in \Gamma\}$ 
is a generating set of $e({\sf riAut}(k))$.  Hence 
$\Gamma = \eta(\Delta)$ for some generating set $\Delta$ of 
${\sf riAut}(k)$.

Moreover, $e(\Delta) = \rho(\tilde{\Gamma}) = $
$\{ (\gamma, {\sf Dom}(\gamma)) : \gamma \in \Gamma\}$, so for every
for every $\delta \in \Delta$,  $e(\delta)$ is of the form 
$(\gamma_{\delta}, {\sf Dom}(\gamma_{\delta}))$ for some 
$\gamma_{\delta} \in \Gamma$; so,
${\sf Dom}(\delta) = {\sf Dom}(\gamma_{\delta})$.  So for
every $\delta \in \Delta$ there exists $\gamma_{\delta} \in \Gamma$ 
such that ${\sf Dom}(\delta) = {\sf Dom}(\gamma_{\delta})$. Thus, if 
$\Gamma$ is such that the map $\rho$ in Theorem \ref{GtildaontoriAut} 
is surjective, then there exists a generating set $\Delta$ as required 
by the surjectiveness condition of \ref{GtildaontoriAut}.  

\smallskip

\noindent (3) \ An example is the four-element generating set of 
$G_{2,1}$ given in \cite{CFP} (pp.\ 240-241); let us call this 
generating set $\Gamma_{\rm CFP}$. The elements of $\Gamma_{\rm CFP}$ 
all have domain codes of cardinality 3 or 4. But any generating set of 
${\sf riAut}(2)$ needs to contain an element with domain 
code of cardinality 2, since composition cannot make domain codes smaller. 

It follows that the elements of $\Gamma_{\rm CFP}$ do not have all the 
domain codes of any generating set of ${\sf riAut}(2)$, so 
$\Gamma_{\rm CFP}$ does not satisfy the surjectiveness condition of 
Theorem \ref{GtildaontoriAut}. 

It follows also that if $\Delta$ is a generating set of 
${\sf riAut}(k)$, then $\eta(\Delta) \neq \Gamma_{\rm CFP}$. Indeed, 
$\Delta$ contains elements of table-size 2 (as we just saw), so 
$\eta(\Delta)$ also has elements of table-size $\leq 2$ (since 
application of $\eta$ means taking the maximum essential extension, 
hence the table-size cannot increase).
But $\Gamma_{\rm CFP}$ has no element of table-size $\leq 2$. 

More strongly, let $\theta$ be any surjective homomorphism
$\theta: (G_{2,1}^{\sim {\cal L}})_{_{\Gamma}} \to {\sf riAut}(k)$.
Then $\Delta = \theta(\tilde{\Gamma})$ is a generating set of 
${\sf riAut}(k)$, hence $\eta(\Delta)$ is a generating set of $G_{2,1}$.
This rules out $\Gamma_{\rm CFP}$, since $\eta(\Delta)$ cannot be equal 
to $\Gamma_{\rm CFP}$.  
 \ \ \ \ \ $\Box$

\newpage 

\begin{pro} \label{existenceofGammaF} \

\noindent {\bf (1)} \ For every generating set $\Delta$ of 
${\sf riAut}_{\sf dict}(k)$ there exists a generating set $\Gamma$ 
of $F_{k,1}$ that satisfies the surjectiveness condition of 
Theorem \ref{GtildaontoriAut} (namely, $\eta(\Delta) \subseteq \Gamma$, 
and for every $\delta \in \Delta$ there exists $\gamma \in \Gamma$ 
with ${\sf Dom}(\delta) = {\sf Dom}(\gamma)$).

If $\Delta$ is finite then $\Gamma$ is finite, and
 \ $|\Delta| \leq 2 \cdot |\Gamma|$.

The generating set $\Gamma = F_{k,1}$ satisfies the surjectiveness 
condition.

\smallskip

\noindent {\bf (2)} \ The condition on $\Gamma$ in Theorem
\ref{GtildaontoriAut} are necessary for the surjectiveness, in general.

\smallskip

\noindent {\bf (3)} \ Not every generating set $\Gamma$ of $F_{k,1}$
satisfies the surjectiveness condition.
More strongly, for some generating set $\Gamma$ of $F_{2,1}$ there is
no surjective homomorphism from $(F_{k,1}^{\sim {\cal L}})_{_{\Gamma}}$ 
onto ${\sf riAut}_{\sf dict}(k)$.
\end{pro}
{\bf Proof.} \ For (1) and (2) the proof is the same as for Proposition
\ref{existenceofGammaG}.  

(3) \ An example is the two-element generating set 
$\{B, \ B^{-1}A\}$ of $F_{2,1}$ derived from the generating set 
$\{A,B\}$ given in \cite{CFP} (pp.\ 222 and 224). The elements have 
domain codes of size 4 (for $B$) or 5 (for $B^{-1}A$). But any generating 
set $\Delta$ of ${\sf riAut}_{\sf dict}(\{a_1,a_2\}^*)$ needs to contain 
an element with domain code of size 3, since composition cannot make 
domain codes smaller. 
It follows that the generating set $\{B, \ B^{-1}A\}$ does not contain 
$\eta(\Delta)$ for any generating set $\Delta$ of 
${\sf riAut}(\{a_1,a_2\}^*)$. 
The rest of the proof is as for Proposition \ref{existenceofGammaG}.
    \ \ \ \ \ $\Box$

\begin{cor} \label{GtildatoriAutCor} \

\noindent $\bullet$ \ The suffix expansion $G_{k,1}^{\sim {\cal L}}$
maps onto ${\sf riAut}(k)$ and $F_{k,1}^{\sim {\cal L}}$ maps onto
${\sf riAut}_{\sf dict}(k)$.

\smallskip

\noindent $\bullet$ \ For every finite generating set $\Delta$ of
${\sf riAut}(k)$ (or of ${\sf riAut}_{\sf dict}(k)$) there exists a 
finite generating set $\Gamma$ of $G_{k,1}$ (respectively $F_{k,1}$) 
with \ $|\Gamma| \ \leq \ 2 \cdot |\Delta|$,
such that the $\Gamma$-generated suffix expansion
$(G_{k,1}^{\sim {\cal L}})_{_{\Gamma}}$ maps onto
${\sf riAut}(k)$ (respectively $(F_{k,1}^{\sim {\cal L}})_{_{\Gamma}}$ 
maps onto ${\sf riAut}_{\sf dict}(k)$).

\smallskip

\noindent $\bullet$ \ There also exist finite generating sets $\Gamma$
of $G_{k,1}$ such that $(G_{k,1}^{\sim {\cal L}})_{_{\Gamma}}$ admits no
surjective homomorphism onto ${\sf riAut}(k)$.
Similarly, there exist finite generating sets $\Gamma$ of $F_{k,1}$ such 
that $(F_{k,1}^{\sim {\cal L}})_{_{\Gamma}}$ has no surjective homomorphism
onto ${\sf riAut}_{\sf dict}(k)$.  
\end{cor}
{\bf Proof.} This follows from Theorem \ref{GtildaontoriAut} and
Propositions \ref{existenceofGammaG}, \ref{existenceofGammaF}.
  \ \ \ \ \ $\Box$

\bigskip

\noindent {\bf Remark.} The fact that ${\sf riAut}(k)$ is a
homomorphic image of  $(G_{k,1}^{\sim {\cal L}})_{_{\Gamma}}$ for
some finite generating set $\Gamma$ of $G_{k,1}$
(Theorem \ref{GtildaontoriAut} and Proposition \ref{existenceofGammaG})
implies that ${\sf riAut}(k)$ is finitely generated (and similarly
for ${\sf riAut}_{\sf dict}(k)$, using Proposition
\ref{existenceofGammaF}).
However, at this point this does not provide a new proof that
${\sf riAut}(k)$ and ${\sf riAut}_{\sf dict}(k)$ are finitely
generated, because we used finite generation in the proofs of
Propositions \ref{existenceofGammaG} and \ref{existenceofGammaF}.

\bigskip

\noindent {\bf Remark.} \ The results in Theorem \ref{GtildaontoriAut}
and Propositions \ref{existenceofGammaG}, \ref{existenceofGammaF} show
that different finite generating sets of $G_{k,1}$ or $F_{k,1}$ can have
very different properties, and the characterization of these finite
generating sets is non-trivial.

\section{Miscellaneous }

\subsection{The map 
$(G_{k,1}^{\sim {\cal L}})_{_{\Gamma}} \to {\sf riAut}(k)$ 
is finite-to-one  }

We give a property of the map $\rho:$
$(G_{k,1}^{\sim {\cal L}})_{_{\Gamma}} \to {\sf riAut}(k) $ showing that
$(G_{k,1}^{\sim {\cal L}})_{_{\Gamma}}$ and ${\sf riAut}(k)$ are
very close.

\begin{pro} \label{rhoFinto1} \
Let $\Gamma$ be any generating set of $G_{k,1}$ (possibly infinite).
The map $\rho$ is {\em finite-to-one}, i.e.,
$\rho^{-1}(\varphi)$ is finite for every $\varphi \in {\sf riAut}(k)$.
It follows that
 \ $(F_{k,1}^{\sim {\cal L}})_{_{\Gamma}} \to $
${\sf riAut}_{\sf dict}(k)$ \ is also finite-to-one.
\end{pro}
{\bf Proof.} \ For any $\varphi \in {\sf riAut}(k)$ let
$(g,S) \in (G_{k,1}^{\sim {\cal L}})_{_{\Gamma}}$ be such that
$\rho(g,S) = \varphi$. Then $g = {\sf max}(\varphi)$, so $g$ is
uniquely determined by $\varphi$. Moreover, for every $h \in S$ we have
${\sf Dom}(\varphi) \subseteq {\sf Dom}(h)$, hence,
${\sf domC}(h) \subset {\sf pref}({\sf domC}(\varphi))$.
Therefore there are only finitely many choices for ${\sf domC}(h)$.
Hence, since $|{\sf imC}(h)| = |{\sf domC}(h)|$ and since
there are only finitely many maximal prefix codes of a given cardinality
(over a given alphabet), there are only finitely many choices for
${\sf imC}(h)$.  Finally, since there are only finitely many bijections
${\sf domC}(h) \to {\sf imC}(h)$, there are only finitely many choices
for $h$.    \ \ \ $\Box$

\subsection{The word problem for the suffix expansion of a group,
and for ${\sf riAut}(k)$ }

We will see that the word problem for the suffix expansion of a monoid 
is closely related to the following problem.
\noindent Let $M$ be a monoid and let $\Gamma$ be a generating set of
$M$. The {\it set word problem} of $M$ over $\Gamma$ is specified as 
follows:

{\sf Input:} Two finite subsets $U = \{u_1, \ldots, u_m\}$ and
$V = \{v_1, \ldots, v_n\}$ of $\Gamma^*$.

{\sf Question:} Is $U = V$ when all strings $u_i, v_j$ are 
evaluated in $M$ ?

\begin{lem} \label{redSetWPtoWP} \  
The set word problem of the monoid $M$ with generating set $\Gamma$ can 
be reduced to the word problem of $M$ over $\Gamma$ by a polynomial-time
{\sc and}-of-{\sc or}'s truth-table reduction.  
\end{lem}
{\bf Proof.} \ We have $U = V$ iff $U \subseteq V$ and $V \subseteq U$.
And we have $U \subseteq V$ iff the following boolean formula is true: 

\smallskip

\hspace{1.1in} 
$\bigwedge_{u_i \in U} \bigvee_{v_j \in V} \, (u_i =_{_M} v_j)$.
 
\medskip 

\noindent This formula involves $n \times m$ calls to the word problem 
of $M$.
 \ \ \ $\Box$

\begin{cor} \label{redWPsuffixexptoWPM} \ 
Let $M$ be a monoid with generating set $\Gamma$. The the word problem
of the suffix expansion $(\tilde{M}^{\cal L})_{\Gamma}$ over $\Gamma$ 
can be reduced to the word problem of $M$ over $\Gamma$ by a 
polynomial-time {\sc and}-of-{\sc or}'s truth-table reduction.
\end{cor}
{\bf Proof.} \ Let ``$x = y$?'' be an input for the word problem of  
$(\tilde{M}^{\cal L})_{\Gamma}$, where $x = x_m \ldots x_1$, 
$y = y_n \ldots y_1$, with $x_i, y_j \in \Gamma$. Let us denote the value 
of a string $w \in \Gamma^*$ in $M$ by $(w)_M$.  
Then the value of $x_m \ldots x_1$ in $(\tilde{M}^{\cal L})_{\Gamma}$ is

\smallskip 

$( (x_m \ldots x_1)_M, \ $
$\{(x_m \ldots x_1)_M, \ \ldots \ , \ (x_2x_1)_M, \ (x_1)_M, \ 1\})$ ,

\smallskip 

\noindent and similarly for $y_n \ldots y_1$. Thus, the word problem of
$(\tilde{M}^{\cal L})_{\Gamma}$  reduces to the conjunction of the word 
problem of $M$ and the set word problem of $M$. 
Lemma \ref{redSetWPtoWP} then yields the result.  
 \ \ \ $\Box$

\begin{cor} \label{WPsuffixexpGk1} \ 
The word problem of the suffix expansion 
$(G_{k,1}^{\sim {\cal L}})_{_{\Gamma}}$  of the Thompson-Higman group
$G_{k,1}$ over a finite generating set $\Gamma$ is in {\sf P}.
\end{cor}
{\bf Proof.} \ This follows from Corollary \ref{redWPsuffixexptoWPM} and 
the fact that the word problem of $G_{k,1}$ over a finite generating set 
is in {\sf P} (proved in \cite{BiThomps}, and strengthened to 
co-context-free in \cite{LehnertSchweitzer}).  
 \ \ \ $\Box$

\begin{pro} \label{wpofriAut} \ 
The word problem of ${\sf riAut}(k)$ over any finite generating set is
in {\sf P}.
\end{pro} 
{\bf Proof.} \ Let $\Delta$ be a finite generating set of 
${\sf riAut}(k)$.  Given a string $x_n \ldots x_1 \in \Delta^*$, the 
table for the value of $x_n \ldots x_1$ in ${\sf riAut}(k)$ can be 
computed by simple composition. 
It was proved in \cite{BiThomps} (Theorem 4.1) that this takes 
polynomial time; in fact it also belongs to the parallel complexity
class ${\sf AC}_1$, which is a subclass of {\sf P}.
 \ \ \ $\Box$


\section{Appendix: \ Monoids of right ideal homomorphisms}

Composing two right-ideal homomorphisms of $A^*$ yields again a
right-ideal homomorphism.
By ${\sf riHom}(k)$ we denote the monoid of all right-ideal homomorphisms
between finitely generated right ideals of $A^*$ (where $|A| = k$), with 
function composition as multiplication.

\begin{lem} \label{rightIdealFacts1} \
For every $\varphi \in {\sf riHom}(k)$, the image ${\sf Im}(\varphi)$ is
a finitely generated right ideal, but
there exists $\varphi \in {\sf riHom}(k)$ such that
$\varphi({\sf domC}(\varphi))$ is not a prefix code.
\end{lem}
{\bf Proof.} \ Let $P = {\sf domC}(\varphi)$ (a finite prefix code) and
${\sf Dom}(\varphi) = PA^*$. Then
${\sf Im}(\varphi) = \varphi(PA^*) =  \varphi(P) \ A^*$, hence
${\sf Im}(\varphi)$ is a finitely generated right ideal.

It is easy to find examples where $\varphi(P)$ is not a prefix code. 
E.g., when $P = \{a, b\}$ and $\varphi$ is defined by the table
$\{ (a,a), (b, aa)\}$, then $\varphi(P) = \{a, aa\}$.
 \ \ \ $\Box$

\bigskip

In \cite{BiRL} (section 3.1) is was proved that 
$\varphi({\sf domC}(\varphi))$ is a prefix code iff the partition 
determined by $\varphi$ on ${\sf Dom}(\varphi)$ is a {\it prefix 
congruence} \footnote{A right-congruence $\equiv$ on a right ideal 
$R \subseteq A^*$ is called a {\it prefix congruence} iff there is
a finite prefix code $P$ and a partition $\equiv_P$ on $P$ such that
for all $x_1, x_2 \in R$: \ $x_1 \equiv x_2 \ \Leftrightarrow \ $
$(\exists w \in A^*) (\exists p_1, p_2 \in P)[ p_1 \equiv_P p_2, $
$x_1 = p_1 w$, $p_2 = x_2 w]$.}. 
This inspires the following.

\begin{defn} \label{defn_riHom_pc} \
Within the monoid ${\sf riHom}(k)$ we define the submonoid

\smallskip

\hspace{0.5in}    ${\sf riHom}_{\sf pc}(k) \ = \ $
$\{ \varphi \in {\sf riHom}(k) \ : \ \varphi({\sf domC}(\varphi)) \ $
${\rm is \ a \ prefix \ code} \}$.

\smallskip

\noindent The elements of ${\sf riHom}_{\sf pc}(k)$ are said to be
{\bf prefix code preserving}.
\end{defn}
The subscript ``{\sf pc}'' stands for ``prefix code''.
It is easy to check that ${\sf riHom}_{\sf pc}(k)$ is indeed a monoid.
The reason for calling the elements of ${\sf riHom}_{\sf pc}(k)$ ``prefix
code preserving'' is the following.

\begin{pro} \label{prefcodePres} \  
For every $\varphi \in {\sf riHom}(k)$ we have:
 \ $\varphi({\sf domC}(\varphi))$ is a prefix code \ iff \ for every 
prefix code $P \subset A^*$, $\varphi(P)$ is a prefix code.
\end{pro}
{\bf Proof.} The right-to-left implication is trivial. To prove the
left-to-right implication, let $x_1, x_2 \in {\sf Dom}(\varphi)$ be 
prefix incomparable, but assume by contradiction that
$\varphi(x_2) = \varphi(x_1) \, w$, for some non-empty $w \in A^*$.
Since $x_1, x_2 \in {\sf Dom}(\varphi)$, there are
$p_1, p_2 \in {\sf domC}(\varphi)$ such that $x_1 = p_1u_1$,
$x_2 = p_2u_2$ (for some $u_1, u_2 \in A^*$). Then
$\varphi(x_2) = \varphi(x_1) \, w$ implies
 \ $\varphi(p_2) \, u_2 = \varphi(p_1) \, u_1w$. This implies that 
$\varphi(p_2)$ and $\varphi(p_1)$ are prefix comparable, which 
contradicts the assumption that $\varphi({\sf domC}(\varphi))$ is a 
prefix code. 
  \ \ \ $\Box$

\medskip

\noindent The following further demonstrates the importance of the monoid
${\sf riHom}_{\sf pc}(k)$.
\begin{pro} \
Every $\varphi \in {\sf riHom}(k)$ has an essential restriction to some
element of ${\sf riHom}_{\sf pc}(k)$.
\end{pro}
{\bf Proof.} It is straightforward to restrict $\varphi$ to some element
$\Phi$ whose image code is ${\sf imC}(\Phi) = A^{\ell}$, where
$\ell$ is the length of a longest string in $\varphi({\sf domC}(\varphi))$.
Obviously, $A^{\ell}$ is a prefix code.
 \ \ \ $\Box$

\bigskip

The Thompson-Higman monoid $M_{k,1}$ (introduced in \cite{BiThomMon}) 
is a homomorphic image of ${\sf riHom}(k)$ and of 
${\sf riHom}_{\sf pc}(k)$.
Indeed, an element of $M_{k,1}$ is an equivalence class of elements of
${\sf riHom}(k)$ or of ${\sf riHom}_{\sf pc}(k)$ where two elements
$\varphi_1$ and $\varphi_2$ are considered equivalent iff they can be
obtained from each other by a finite number of essentially equal
restrictions and essentially equal extensions.

\bigskip

As a generalization of the monoid ${\sf riAut}(k)$ that we introduced 
earlier, we consider the monoid ${\sf riIso}(k)$ consisting of
all right-ideal isomorphisms between finitely generated ideals (not
necessarily essential) of $A^*$.  
For ${\sf riIso}(k)$ and in particular, for ${\sf riAut}(k)$, we have:

\begin{pro} \label{pcFORriIso} \
Every element $\varphi \in {\sf riIso}(k)$ is prefix-code preserving.
\end{pro}
{\bf Proof.} \ For $\varphi \in {\sf riIso}(k)$ let
$P = {\sf domC}(\varphi)$ (a finite prefix code) and let
$Q = \varphi(P)$. If $Q$ is not a prefix code then there exist
$q_1 \neq q_2 \in Q$ with $q_2 = q_1v$ for some $v \in A^*$,
$v \neq \varepsilon$. Since $\varphi$ is an injective homomorphism there
exist $p-1 \neq p-2 \in P$ such that
$q_1 = \varphi(p_1) \neq \varphi(p_2) = q_2 = \varphi(p_1) \ v$
$ =  \varphi(p_1 v)$. By injectiveness, $p_2 = p_1v$, which contradicts
the fact that $P$ is a prefix code.
 \ \ \ $\Box$

\bigskip

As a consequence, ${\sf riIso}(k)$ consists of all right-ideal
isomorphisms $\varphi$ such that ${\sf domC}(\varphi)$ and
$\varphi({\sf domC}(\varphi)) = {\sf imC}(\varphi)$ are prefix codes
(not necessarily maximal).
And ${\sf riAut}(k)$ consists of all right-ideal isomorphisms $\varphi$
such that ${\sf domC}(\varphi)$ and
$\varphi({\sf domC}(\varphi)) = {\sf imC}(\varphi)$ are maximal prefix
codes. The notation {\sf riAut}, where ``Aut'' stands for automorphism,
is motivated by the fact that ${\sf riAut}(k)$ maps onto the group
$G_{k,1}$. 

Note that ${\sf riIso}(k)$ does not map onto $G_{k,1}$. Indeed, 
${\sf riIso}(k)$ has a zero (the empty map), the only group that 
${\sf riIso}(k)$ maps onto is the one-element group.
 

\bigskip

\bigskip



\bigskip

\bigskip

\noindent {\bf Jean-Camille Birget} \\
Dept.\ of Computer Science \\
Rutgers University at Camden \\
Camden, NJ 08102, USA \\
{\tt birget@camden.rutgers.edu}

\end{document}